\documentclass[11pt,psamsfonts]{amsart}

\usepackage{a4wide}
\usepackage{amssymb}
\usepackage{amstext}
\usepackage[ansinew]{inputenc}

\let\ge=\geqslant
\let\geq=\geqslant
\let\le=\leqslant
\let\leq=\leqslant
\let\ptl=\partial
\let\Sg=\Sigma
\let\sg=\sigma
\let\eps=\varepsilon
\let\Om=\Omega

\newcommand{\rr}{\mathbb{R}}
\newcommand{\hh}{\mathbb{H}}
\newcommand{\h}{\mathcal{H}}
\newcommand{\pp}{\mathcal{P}}
\newcommand{\nn}{\mathbb{N}}

\newcommand{\wt}{\widetilde}

\DeclareMathOperator{\vol}{vol}

\DeclareMathOperator{\ric}{Ric}

\newtheorem{theorem}{Theorem}[section]
\newtheorem{proposition}[theorem]{Proposition}
\newtheorem{lemma}[theorem]{Lemma}
\newtheorem{corollary}[theorem]{Corollary}

\theoremstyle{definition}

\newtheorem{remark}[theorem]{Remark}

\theoremstyle{remark}

\newenvironment{enum}{\begin{enumerate}
}{\end{enumerate}}

\begin{document}

\title[Isoperimetric comparison theorems for convex bodies]{Some isoperimetric 
comparison theorems for convex bodies in Riemannian manifolds}

\author[V.~Bayle]{Vincent Bayle}
\address{Institut Fourier \\BP 74 \\38402 Saint Martin D'heres Cedex \\ France}
\email{vbayle@ujf-grenoble.fr}

\author[C.~Rosales]{C\'esar Rosales}
\address{Departamento de Geometr\'{\i}a y Topolog\'{\i}a \\
Universidad de Granada \\ E--18071 Granada \\ Espa\~na}
\email{crosales@ugr.es}

\date{November 7, 2003}
\subjclass{53C20, 49Q20}
\keywords{Isoperimetric profile, isoperimetric regions, differential inequality, 
comparison theorems}

\begin{abstract} 
We prove that the isoperimetric profile of a convex domain $\Om$ with 
compact closure in a Riemannian manifold $(M^{n+1},g)$ satisfies a 
second order differential inequality which only depends on the dimension 
of the manifold and on a lower bound on the Ricci curvature of $\Omega$. 
Regularity properties of the profile and topological 
consequences on isoperimetric regions arise naturally from this differential 
point of view. 

Moreover, by integrating the differential inequality we 
obtain sharp comparison theorems: not only can we derive an inequality 
which should be compared with L\'evy-Gromov Inequality but we also show 
that if $\ric\geq n\delta$ on $\Om$, then the profile of $\Om$ is bounded 
from above by the profile of the half-space $\mathbb{H}_{\delta}^{n+1}$ in 
the simply connected space form with constant sectional curvature $\delta$. 
As consequence of isoperimetric comparisons we obtain geometric 
estimations for the volume and the diameter of $\Om$, and for the first non-zero 
Neumann eigenvalue for the Laplace operator on $\Om$. 
\end{abstract}

\maketitle

\thispagestyle{empty}

\section{Introduction}
\label{sec:introduction}

Let $\Om$ be a domain (connected open set) with non-empty boundary of 
a Riemannian manifold $(M^{n+1},g)$. The so-called \emph{partitioning problem} 
in $\Om$ consists on finding, for a given $V<\vol(\Om)$, a 
minimum of the perimeter functional $\pp(\cdot\, ,\Om)$ in the class 
of sets in $\Om$ that enclose volu\-me $V$. Here $\vol(E)$ is the 
$(n+1)$-dimensional Hausdorff measure of a set $E\subseteq M$ and 
$\pp(E,\Om)$ denotes the \emph{perimeter} of $E$ \emph{relative to} 
$\Om$, which essentially measures the area of $\ptl E\cap\Om$ (see
Section \ref{sec:preliminaries} for a precise definition). Solutions  
to the partitioning problem are called \emph{isoperimetric regions} or 
\emph{minimizers} in $\Om$ of volume $V$.  

The partitioning problem is object of an intensive study. The first 
questions taken into consideration were related to the existence 
and regularity of minimizers. In the light of standard results in 
Geometric Measure Theory \cite{morgan}, inside a smooth domain $\Om$ 
with compact closure, minimizers do exist for any given volume and 
their boundaries are smooth, up to a closed set of singularities with 
high Hausdorff codimension, (see Proposition \ref{prop:regularity} 
for a precise statement). Recently, geometric and topological properties 
of minimizers have been studied by A.~Ros and E.~Vergasta~\cite{rosvergasta} 
and P.~Sternberg and K.~Zumbrun~\cite{zumbrun2} inside a Euclidean convex
body, and by M.~Ritor\'e and C.~Rosales \cite{ritros} inside Euclidean cones. 
However, in spite of the last advances, the complete description 
of isoperimetric regions has been achieved only for certain convex 
domains such as half-spaces in the simply connected space forms, 
Euclidean balls, Euclidean slabs, and Euclidean convex cones, among 
others. A beautiful survey containing most of the results above, 
including recent progress and open questions is the one by 
A.~Ros \cite{ros}. 

Much of the information concerning the partitioning problem is 
contained in the \emph{isoperimetric profile of} $\Om$: the function 
$I_\Om(V)$ which assigns to $V$ the least-perimeter separation of 
volume $V$ in $\Om$. In this paper, assuming that $\Om$ is a convex 
domain with compact closure in $(M^{n+1},g)$, we prove regularity 
properties of the profile, connectivity results for minimizers and for 
their boundaries, and above all, we obtain sharp lower and upper bounds 
for the isoperimetric profile involving the infimum of the Ricci 
curvature of $\Om$.

We begin this work with a preliminary section where we introduce 
the notation and give some basic results. For example, Proposition 
\ref{asbpnz} is an adaptation to the partitioning problem of a 
result by P.~B\'{e}rard and D.~Meyer \cite{bemeyer} in which it is 
shown that the isoperimetric profile $I_\Om$ approaches asymptotically 
the profile of the half-space in $\rr^{n+1}$ for small volumes. We also 
summarize existence and regularity results for isoperimetric regions in 
Proposition~\ref{prop:regularity}, and state an analytic comparison 
result for the solutions of a differential inequality (Theorem 
\ref{th:analcomp}) that will be useful in Section~\ref{compmod}. 

In Section \ref{sec:results}, inspired by previous results by C.~Bavard 
and P.~Pansu \cite{BP}, P.~Sternberg and K.~Zumbrun \cite{zumbrun2}, 
F.~Morgan and D.~Johnson \cite{mj}, and V.~Bayle \cite{Ba}, we prove 
(Theorem \ref{prop:inequality}) that the renormalized isoperimetric 
profile $Y_{\Omega}=I_{\Omega}^{(n+1)/n}$ of a smooth convex domain 
$\Om$ with compact closure satisfies a second order differential 
inequality of the type
\begin{equation}
\label{eq:bombero}
Y_\Om''\leq C\,Y_\Om^{(1-n)/(1+n)},
\end{equation}
where $C$ is a constant depending on the dimension of the ambient manifold 
and on a lower bound on the Ricci curvature over $\Om$.  

The idea of the proof of \eqref{eq:bombero} relies on a local comparison 
of $Y_\Om$ with the renormalized profile $P(V)^{(n+1)/n}$ associated 
to the deformation of a minimizer $E$ -which exists by the compactness of 
$\overline{\Om}$- given by equidistant hypersurfaces 
to $\ptl E\cap\Om$. Some technical difficulties arise due 
to the possible presen\-ce in high dimensions of singularities in 
$\ptl E\cap\Om$. These difficulties are solved by an approximation 
argument consisting in the construction of ``almost parallel variations" 
(Lemma~\ref{lem:phieps}). 
This scheme of proof was previously used in \cite{mj} and \cite{Ba} to get 
a differential inequality for the isoperimetric profile of a closed 
Riemannian manifold, and in \cite{ritros} to characterize isoperimetric 
regions in smooth convex cones. As in \cite{ritros}, our proof differs from those 
of \cite{mj} and \cite{Ba} in the presence of a boundary term involving the 
second fundamental form of $\ptl\Om$ which can be controlled by using the 
convexity of $\Om$. 

From the differential inequality \eqref{eq:bombero}, that yields concavity 
of the profile under the assumption of non-negative Ricci curvature on $\Om$ 
(Theorem~\ref{thconpc}), we derive regularity properties of the profile 
(Proposition~\ref{syprdn}) and topological consequences related to the 
connectivity of minimizers and isoperimetric hypersurfaces 
(Propositions \ref{prop:syprdn} and \ref{syprdn2}). Similar 
previous results for closed Riemannian manifolds and for convex bodies in
the Euclidean setting were established in \cite{BP}, \cite{mj}, 
\cite{zumbrun2}, \cite{K} and \cite{Ba}. 

In Section \ref{compmod} we use analytic arguments to
obtain geometric comparison theorems. As a matter of fact, 
integration of the differential inequality \eqref{eq:bombero} 
makes possible to compare the isoperimetric profile of a smooth 
convex domain $\Om$ with compact closure and $\ric\geq n\delta$, 
with an exact solution of the differential equation associated to 
\eqref{eq:bombero} that satisfies either the same initial conditions 
or the same boundary conditions. On the one hand we prove 
in Theorem \ref{thcomp} the isoperimetric inequality
\begin{equation}
\label{eq:bombero2}
I_\Om(V)\leq I_{\hh_\delta^{n+1}}(V),\qquad V\in [0,\vol(\Om)],
\end{equation}
where $\mathbb{H}^{n+1}_\delta$ is a half-space in the simply connected 
space form with constant sectional curvature $\delta$. In Remark~\ref{re:mj} 
we indicate that the geometric arguments employed by F.~Morgan and 
D.~Johnson in \cite[Theorem 3.5]{mj} can be adapted to prove that 
\eqref{eq:bombero2} is also valid for unbounded convex domains. In 
Theorem~\ref{th:unbounded} we show that equality in \eqref{eq:bombero2} for 
some $V_0\in (0,\vol(\Om)]$ implies that $\ptl\Om$ is a totally geodesic 
hypersurface and $\Om$ has constant sectional curvature $\delta$ in a 
neighborhood of $\ptl\Om$. 

On the other hand, in Corollary \ref{cgenerallevygromov} we deduce a lower
bound for the profile that should be compared with L\'evy--Gromov 
inequality \cite{GR2}. In precise terms, we prove that any Borel set $E$ 
contained in a smooth convex body $\Om$ with $\ric\geq n\delta>0$, satisfies
\begin{eqnarray}
\label{eq:levygromov2}
\frac{\mathcal{P}(E,\Omega)}{\vol(\Omega)}\geq 
\frac{\mathcal{P}(E^*,\mathbb{H}^{n+1}_{\delta})}
{\vol(\mathbb{H}^{n+1}_{\delta})},
\end{eqnarray}
where $E^*\subseteq\mathbb{H}^{n+1}_\delta$ is a half-ball centered 
at $\ptl\mathbb{H}^{n+1}_{\delta}$ with
$\vol(E)/\vol(\Omega)=\vol(E^*)/\vol(\mathbb{H}^{n+1}_{\delta})$. 
Moreover, inequality \eqref{eq:levygromov2} is sharp since equality 
for a proper set $E\subset\Om$ implies that $\Om$ is isometric to 
$\hh_\delta^{n+1}$. 

Our isoperimetric inequalities in Section \ref{compmod} can be used, 
as in \cite{G1} and \cite{Ba}, to derive comparison theorems for convex 
bodies involving geometric quantities such as the volume or the 
diameter, see Theorem~\ref{th:bishop}, Remark~\ref{re:altproofca} and 
Theorem~\ref{th:cheng}. Furthermore, by reproducing the symmetrization 
arguments in \cite[Th\'eor\`eme 5]{bemeyer} we prove in 
Theorem~\ref{theocompeigenvalue} that if $\ric\geq n\delta>0$ on $\Om$, 
then the lowest non-zero eigenvalue for the Laplace operator in $\Om$ 
with Neumann boundary condition is bounded from below by the one of the
half-sphere $\hh_{\delta}^{n+1}$ of radius $1/\sqrt{\delta}$, with equality 
if and only if $\Om$ is isometric to~$\hh_\delta^{n+1}$.   

Finally, we have added in a last section as an appendix a geometric proof 
of inequality \eqref{eq:bombero2} for the case of a smooth convex body 
$\Om$ in $\rr^{n+1}$.

As mentioned in \cite{zumbrun2}, in addition to the geometric interest 
of this work, we remark that the partitioning problem can be linked 
with a well-studied variational question related to phase transitions 
(see also \cite{zumbrun2bis}). 

\vspace{0,1cm}
\noindent
\emph{Acknowledgements.} The idea of this work was conceived while 
V.~Bayle was visiting the~ University of Granada in the spring of 2003. 
The first author was supported by the Marie~Curie Research Training 
Networks ``EDGE", HPRN-CT-2000-00101. The second author was supported 
by MCyT-Feder research project BFM2001-3489. Both authors express their 
deep thanks to Manuel Ritor\'{e} for his encoura\-gement and helpful 
comments during the preparation of these notes.

\section{Preliminaries}
\label{sec:preliminaries}
\setcounter{equation}{0}

\subsection{The isoperimetric profile}
Let $\Om$ be a smooth domain (connected open set) with compact closure 
$\overline{\Om}$ contained in a Riemannian manifold $(M^{n+1},g)$. The 
$(n+1)$-dimensional and the $k$-dimensional Hausdorff measures
of a Borel set $E\subseteq M$ will be denoted by $\vol(E)$ and $\h_k(E)$ 
respectively. For any measurable set $E\subseteq M$, let $\pp(E,\Om)$ be 
the De Giorgi~\emph{perimeter} 
of $E$ \emph{relative to} $\Om$, defined as
\[
\pp(E,\Om)=\sup\left\{\int_E\text{div}\,Y\,d\h_{n+1}:g(Y,Y)\leq1\right\},
\] 
where $Y$ is a smooth vector field over $M$ with compact support contained 
in $\Om$, and $\text{div}\,Y$ is the divergence of $Y$ \cite[p. 140]{chavel}. 
If, for instance, $E$ has $C^2$ boundary, then $\pp(E,\Om)=\h_n(\ptl E\cap\Om)$ 
by the Gauss--Green theorem. 

A set $E\subseteq M$ is said to be of \emph{finite perimeter in} 
$\Om$ if $\pp(E,\Om)<\infty$. We refer to the reader to \cite{giusti}, 
\cite{bvfunctions} and \cite{C2} for background about perimeter, sets of 
finite perimeter, and their use in the context of the Geometric Measure 
Theory.

The \emph{isoperimetric profile} of $\Om$ is the function $I_\Om:[0,\vol(\Om)]\to\rr^{+}\cup\{0\}$ given by 
\[
I_\Om(V)=\inf\,\{\pp(E,\Om):\,E\subseteq\Om,\,\, \vol(E)=V\},
\]
where the infimum is taken over sets of finite perimeter in $\Om$. 
We define the \emph{renormalized isoperimetric profile} of $\Om$ as 
the function $Y_\Om=I_\Om^{(n+1)/n}$.

\vspace{0,1cm}

Through this paper we shall use the following basic properties of the 
isoperimetric profile
\begin{itemize}
\item[$\bullet$] $I_\Om$ is a non-negative function which only vanishes at 
$V=0$ and $V=\vol(\Om)$.
\vspace{0,05cm}
\item[$\bullet$] $I_\Om(V)=I_\Om(\vol(\Om)-V)$, \quad $V\in [0, \vol(\Om)]$. 
\item[$\bullet$] $I_\Om$ is a lower semicontinuous function 
\cite[Theorems 1.9 and 1.19]{giusti}.
\end{itemize}

\vspace{0,1cm}

The following proposition is an adaptation of a result by 
P.~B\'{e}rard and D.~Meyer \cite[App. C]{bemeyer}, in which 
the cited authors show that the isoperimetric profile of a 
closed manifold $(M^{n+1},g)$ (i.e., a compact Riemannian 
manifold without boundary) asymptotically approaches the 
profile of $\rr^{n+1}$ for small volumes. 

\begin{proposition} 
\label{asbpnz}
Let $\Om$ be a smooth domain with compact closure and non-empty 
boundary in a Riemannian manifold $(M^{n+1},g)$. Denote by 
$\mathbb{H}^{n+1}$ the half-space $\{x_{n+1}>0\}$ in $\rr^{n+1}$. 
Then, the asymptotic behaviour of the isoperimetric profile of 
$\Om$ at the origin is 
\begin{eqnarray*}
I_{\Omega}(V)\underset{\substack{V\to 0 \\ V>0 }}
{\sim}I_{\mathbb{H}^{n+1}}(V)= 2^{-1/(n+1)}\,\gamma_{n+1}\, V^{n/(n+1)},
\end{eqnarray*}
where $\gamma_{n+1}=\h_n(\mathbb{S}^n)/\h_{n+1}(B(1))^{{n}/(n+1)}$ 
stands for the $(n+1)$-dimensional Euclidean isoperimetric constant.

As a consequence, the right derivative of the renormalized profile 
at the origin is given by
\[
\label{eq:righder}
(Y_{\Omega})'_r(0)=2^{-1/n}\,\gamma_{n+1}^{(n+1)/n}.
\]
\end{proposition}

\begin{proof}
The only change with respect to the proof by P.~B\'{e}rard and D.~Meyer 
that must be taken into account consists in proving a localization 
lemma as in \cite[p. 531]{bemeyer} for any small geodesic ball $B$ 
centered at $\ptl\Om$ and intersected with $\Om$. In precise terms, we 
need to show that inside $B\cap\Om$ the isoperimetric inequality 
for the relative perimeter infinitesimally behaves as in 
$\mathbb{H}^{n+1}$. This property comes from the fact that $B\cap\Om$ 
is almost isometric to a half-ball in $\mathbb{H}^{n+1}$ centered at 
$\ptl\mathbb{H}^{n+1}$. Finally, a compactness argument as in 
\cite{bemeyer} allows us to pass from the localization lemmae to a 
global isoperimetric inequality. 
\end{proof}

\begin{remark}
\label{re:conito} 
The asymptotic behaviour in the proposition above provides upper and 
lower bounds on the profile for small volumes. In fact, for any 
$\eps>0$, there exists $V(\Omega,\eps)>0$ such that 
\[
(1-\eps)\,I_{\mathbb{H}^{n+1}}(V)\leq I_\Om(V)\leq 
(1+\eps)\,I_{\mathbb{H}^{n+1}}(V),\quad\text{ whenever } V\leq 
(\Omega,\eps).
\]
The last inequality and the one given in \cite[App. C]{bemeyer} imply 
that a set $E$ in $\Om$ such that $\pp(E,\Om)=I_\Om(V)$ for a small 
volume $V$, must touch the boundary of $\Om$.  
\end{remark}

Now, we introduce another notion of isoperimetric profile 
(see \cite{GR2}, \cite{G1} and \cite{Ba}), which is sometimes 
more relevant in order to obtain comparison theorems. It is 
given by the function $h_{\Omega}:[0,1]\to\rr^{+}\cup\{0\}$, 
defined for all $\beta$ in $[0,1]$ by
\begin{equation}
\label{rmqopr}
h_{\Omega}(\beta)=\frac{I_{\Omega}\big(\beta\vol(\Omega)\big)}{\vol(\Omega)}.
\end{equation}
This point of view, which somehow corresponds to the choice of 
a probability measure on $\Omega$, will be considered in the 
proof of a L\'evy-Gromov type inequality 
(Theorem \ref{generallevygromov}).

\subsection{Isoperimetric regions: existence and regularity}
\label{subsec:eandr}

Let $\Om$ be a smooth domain of a Riemannian manifold $(M^{n+1},g)$. An 
\emph{isoperimetric region} -or simply a \emph{minimizer}- in $\Om$ for 
volume $V\in (0,\vol(\Om))$ is a set $E\subseteq\Om$ such that 
$\vol(E)=V$~and $\pp(E,\Om)=I_\Om(V)$. 

In the following proposition we summarize some results from Geometric 
Measure Theory concerning the existence and regularity of isoperimetric 
regions in $\Om$. 

\begin{proposition}[\cite{giusti}, \cite{gmt}, \cite{gruter}, \cite{morgan2}, 
\cite{bombieri}]
\label{prop:regularity}
Let $\Om$ be a smooth domain with compact closure in a Riemannian manifold 
$(M^{n+1},g)$. For any $V\in (0,\vol(\Om))$ there is an open set 
$E\subset\Om$ which minimizes the perimeter relative to $\Om$ for volume 
$V$. The boundary $\Lambda=\overline{\ptl E\cap\Om}$ can be written as a 
disjoint union $\Sg\cup\Sg_0$, where $\Sg$ is the regular part of $\Lambda$ 
and $\Sg_0=\Lambda-\Sg$ is the set of~singularities. Precisely, we have
\begin{enum}
\item $\Sg\cap\Om$ is a smooth, embedded hypersurface with constant
mean curvature.  
\item If $p\in\Sg\cap\ptl\Om$, then $\Sg$ is a smooth, embedded hypersurface
with boundary contained in $\ptl\Om$ in a neighborhood of $p$; in this
neighborhood $\Sg$ has constant mean curvature and meets $\ptl\Om$
orthogonally.
\item $\Sg_0$ is a closed set of Hausdorff dimension less than or
equal to $n-7$.
\item At every point $q\in\Sg_0$ there is a tangent minimal cone 
$C\subset T_qM$ different from a hyperplane. The square sum 
$|\sg|^2=k_1^2+ \ldots+k_n^2$ of the principal curvatures of $\Sg$ 
tends to $\infty$ when we approach $q$ from $\Sg$.  
 \end{enum}
\end{proposition}

In the preceding proposition the regular set $\Sg$ is defined as
follows: for $p\in\Sg$ there is a neighborhood $W$ of $p$ in $\Sg$
such that $W$ is a smooth, embedded hypersurface without boundary or
with boundary contained in $\ptl\Om$. Note that a consequence
of the proposition above is the absence of interior points in $\Sg$
meeting $\ptl\Om$ tangentially, see \cite{gruter-nuevo}.

\begin{remark}
\label{re:voidboundary}
The regular hypersurface $\Sg$ associated to a minimizer in $\Om$ 
need not meet the boundary of $\Om$. An example illustrating this 
situation can be found at the end of Section 2 in \cite{ritros}.
\end{remark}

\subsection{An analytic comparison result}
\label{subsec:analytic}

Let $f:I\to\rr$ be a function defined on an open interval. For any 
$x_0\in I$ we denote by $\overline{D^2f}(x_0)$ the 
\emph{upper second derivative} of $f$ at $x_0$, defined by 
\begin{equation}
\label{eq:uppersecondderivative}
\overline{D^2f}(x_0)=\limsup_{h\to 0}\,
\frac{f(x_0+h)+f(x_0-h)-2f(x_0)}{h^2}.
\end{equation}

The main tool that we shall employ in Section \ref{compmod} to derive 
comparison theorems from differential inequalities is the following 
technical result. A detailed development is included in \cite{Bat} and 
will appear in a forthcoming paper.

\begin{theorem}
\label{th:analcomp}
Let $f,g:[0,a]\to \rr$ be continuous functions with positive values on 
$(0,a)$. Let $H:\rr^+\to\rr$ be the function 
$H(x)=-\alpha\delta x^{(2-\alpha)/\alpha}$, where $\delta\in\rr$ and 
$\alpha\geq 2$. Suppose that $f$ satisfies the second order differential 
inequality
\[
\overline{D^2f}(x)\leq H[f(x)],\qquad x\in (0,a),
\]
while $g$ is a $C^2$-function that satisfies the differential equation
\[
g''(x)=H[g(x)],\qquad x\in (0,a).
\]
Then, we have
\begin{itemize}
\item[(i)] If $f(0)=g(0)$ and $f(a)=g(a)$, then $f\geq g$ on $[0,a]$. 
Moreover, if $f(x_0)=g(x_0)$ for some $x_0\in (0,a)$, then $f=g$ on 
$[0,a]$.
\item[(ii)] If $f(0)=g(0)$ and the right derivatives at the origin 
verify $f'_r(0)\leq g'_r(0)<+\infty$, then $f\leq g$ on $[0,a]$. 
Moreover, if $f(x_0)=g(x_0)$ for some $x_0\in (0,a]$, then $f=g$ 
on $[0,x_0]$.
\end{itemize}
\end{theorem} 

The theorem above can be seen as a generalization of the fact that a 
concave function $f$ on $[0,a]$ is pinched between any tangent line 
and the secant line passing through $(0,f(0))$ and $(a,f(a))$.

\subsection{Convex domains in Riemannian manifolds}
The term ``convex domain" is used in different non-equivalent ways 
in the literature. We adopt the following definition:

Let $\Om$ be a domain of a Riemannian manifold $(M^{n+1},g)$. We say 
that $\Om$ is \emph{convex} if any two points 
$p,q\in\Om$ can be joined by a minimizing geodesic of $M$ which is 
contained in $\Om$. A convex domain $\Om$ with compact closure in $M$ 
will be called a \emph{convex~body}. 

The convexity of a smooth domain $\Om$ implies the \emph{local convexity} 
of $\ptl\Om$, which means that all the geodesics in $M$ tangent to 
$\ptl\Om$ are locally outside of $\Om$. As R.~Bishop proved (\cite{convexity}),
the local convexity of $\ptl\Om$ is equivalent to an analytic condition 
(the so-called \emph{infinitesimal convexity}) involving the second 
fundamental form of $\ptl\Om$. As a consequence, we have
\begin{quotation}
``\emph{A smooth convex domain $\Om$ of a Riemannian manifold 
$\emph{(}M^{n+1},g\emph{)}$ satisfies that the second 
fundamental form $\emph{II}_p$ of $\ptl\Om$ with respect to the 
normal pointing into $\Om$ is positive semidefinite at any 
$p\in\ptl\Om$}''.
\end{quotation}

\begin{remark}
Most of the results of the paper in which the convexity of $\Om$ 
is assumed are also valid under the weaker condition that $\text{II}_p$ 
is positive semidefinite at any $p\in~\!\!\ptl\Om$.
\end{remark}

The following result is an application to the setting of convex 
bodies of two well-known comparison theorems in Riemannian 
Geometry. It will be useful in order to show that our isoperimetric 
inequalities in Section \ref{compmod} are sharp.

\begin{theorem}
\label{th:bishop}
Let $\Om$ be a smooth convex domain of a complete Riemannian 
manifold~$(M^{n+1},\newline g)$. Denote by $\hh^{n+1}_\delta$ 
the $(n+1)$-dimensional half-sphere of radius $1/\sqrt{\delta}$. 
If the Ricci curvature of $M$ satisfies $\ric\geq n\delta>0$ 
on $\Om$, then
\begin{itemize}
\item[(i)] $\overline{\Om}$ is compact and 
$\emph{diam}(\Om)\leq\pi/\sqrt{\delta}$ \emph{(}Bonnet--Myers 
Theorem\emph{)}.
\item[(ii)] If $\ptl\Om\neq\emptyset$ then $\vol(\Om)\leq
\vol(\hh^{n+1}_\delta)$ and equality implies that $\ptl\Om$ 
is totally geodesic in $M$ and $\Om$ is isometric to 
$\hh^{n+1}_\delta$ \emph{(}Bishop's Theorem\emph{)}. 
\end{itemize}
\end{theorem}

\begin{proof}
Assertion (i) is a consequence of the Bonnet--Myers theorem 
(\cite[Theorem 2.12]{chavel}). Now we give an outline of the 
proof of (ii) which follows the arguments in \cite[Theorem 3.9]{chavel}. 
Fix $p\in\ptl\Om$ and call $\nu$ to the unit normal vector to $\ptl\Om$ 
at $p$ pointing into $\Om$. Denote 
$\nu^+=\{\eta\in T_p M: g_p(\eta,\eta)=1,\,\,g_p(\eta,\nu)>0\}$. For any 
$\eta\in\nu^+$ let $\gamma_\eta$ be the unique geodesic in $M$ such that 
$\gamma_\eta(0)=p$ and $\gamma_\eta'(0)=\eta$. Let also denote
\begin{align*}
c(\eta)&=\sup\,\{t>0: \text{dist}(p,\gamma_\eta(t))=t\}, 
\\
d(\eta)&=\sup\,\{t>0: \gamma_\eta((0,t])\subset\Om\}.
\end{align*}
Let $C(p)$ be the cut locus of $p$ in $M$. By the convexity of $\Om$ we 
have
\[
\Om\setminus C(p)=\{\gamma_\eta(t):\eta\in\nu^+,\,\,t\in (0,m(\eta))\},
\] 
where $m(\eta)=\min\{c(\eta),d(\eta)\}$ is continuous as function of 
$\eta\in\nu^+$. Call $J(t,\eta)$ to the~Jacobian determinant of the 
map $(t,\eta)\in (0,\infty)\times\nu^+\mapsto\gamma_\eta(t)$. By using 
the integration formula for polar geodesic coordinates around $p$ and 
Bishop's Theorem \cite[Theorem 3.8]{chavel}, we get
\begin{align}
\label{eq:perto}
\vol(\Om)&=\int_{\nu^+}\left (\int_0^{m(\eta)}J(t,\eta)\,dt\right)\,
d\eta\leq\int_{\nu^+}\left(\int_0^{m(\eta)}S_\delta(t)^n\,dt\right)\,d\eta 
\\
\nonumber
&\leq\int_{\nu^+}\left(\int_0^{\pi/\sqrt{\delta}}S_\delta(t)^n\,dt\right)\,
d\eta=
\vol(\mathbb{H}^{n+1}_\delta),
\end{align}
where $S_\delta(t)=\sin(\sqrt{\delta}\,t)/\sqrt{\delta}$ and we have 
used (i) to write $m(\eta)\leq\pi/\sqrt{\delta}$. 

If equality holds in \eqref{eq:perto}, then $m(\eta)=\pi/\sqrt{\delta}$ 
for all $\eta\in\nu^+$ and $\Om\setminus C(p)$ has constant 
sectional curvature $\delta$ with respect to any plane containing a 
tangent vector to a geodesic $\gamma_\eta$. As consequence, 
$\Om\setminus C(p)$ is a geodesic half-ball $B$ centered at $p$ with 
radius $\pi/\sqrt{\delta}$ and we deduce, as in the local Cartan's 
theorem \cite[Exercise 3.1]{chavel}, that $B$ is isometric to 
$\mathbb{H}^{n+1}_\delta$. From this, it is easy to see that 
$\Om\cap C(p)=\emptyset$ and $\ptl\Om=\{\gamma_\eta(t):\eta\in 
T_p(\ptl\Om),\,\,t\in [0,\pi/\sqrt{\delta}]\}$.
\end{proof}

\begin{remark}
In the theorem above we assume $\ric\geq n\delta>0$ only in $\Om$. 
In Section \ref{compmod} we shall see that the two geometric inequalities 
in Theorem~\ref{th:bishop} can be obtained by using isoperimetric 
comparisons. In Theorem~\ref{th:cheng} we characterize the half-spheres 
as the only convex domains for which equality in 
Theorem~\ref{th:bishop} (i) holds.  
\end{remark}

\section{The differential inequality}
\label{sec:results}
\setcounter{equation}{0}

Let $\Om$ be a smooth convex body of a Riemannian manifold $(M^{n+1},g)$. Our main goal in 
this section is to prove that the renormalized isoperimetric profile $Y_\Om$ satisfies a 
differential inequality as that as in \eqref{eq:bombero}. We shall then derive some immediate consequences related to the regularity of the profile and the connectivity of isoperimetric regions in $\Om$.

Let us start with the proof of the differential inequality. As we pointed out in Section 
\ref{sec:introduction}, the idea of the proof consists in a local comparison of $Y_\Om$ 
with the relative profiles associated to ``almost parallel variations" of a minimizer $E$ 
in $\Om$ for a fixed volume $V_0$. These variations will be constructed by using the following lemma

\begin{lemma}
\label{lem:phieps}
Let $E$ be an isoperimetric region inside a smooth domain $\Om$ with compact closure in a Riema\-nnian manifold $(M^{n+1},g)$. Denote by $\Sg$ the regular part of $\Lambda=\overline{\ptl E\cap\Om}$. Then, there is a sequence $\{\varphi_\eps:\Sg\to\rr\}_{\eps>0}$ of smooth functions with compact 
support in $\Sg$, such that 
\begin{itemize}
\item[(i)] $0\leq\varphi_\eps\leq 1$,\quad $\eps>0$.
\item[(ii)] $\{\varphi_\eps\}\to 1$ in the Sobolev space $H^1(\Sg)$, that is
\[
\lim_{\eps\to 0}\,\int_\Sg\varphi_\eps^2\,d\h_n=\pp(E,\Om),\qquad
\lim_{\eps\to 0}\,\int_\Sg|\nabla\varphi_\eps|^2\,d\h_n=0,
\]
where $\nabla\varphi_\eps$ is the gradient of $\varphi_\eps$ relative to $\Sg$.
\item[(iii)] $\lim_{\eps\to 0}\,\varphi_\eps(p)=1$,\quad $p\in\Sg$.
\end{itemize}  
\end{lemma} 

A complete proof of the lemma above when $\Om$ is a Euclidean domain 
can be found in \cite[Lemma 2.4]{zumbrun2}. The general case is treated 
in a similar way, see \cite{morganritore} and  
\cite[Proposition 1.1]{Ba} for further details. In 
\cite[Lemma 3.1]{morganritore} it was shown that the existence of 
$\{\varphi_\eps\}_{\eps>0}$ is guaranteed for a bounded, constant 
mean curvature hypersurface $\Sg$ with a closed singular set 
$\Sg_0=\overline{\Sg}-\Sg$ such that $\h_{n-2}(\Sg_0)=0$ or consisting 
of isolated points. 

Now, we can prove the main result of this section. Recall that 
$\overline{D^2f}(x_0)$ denotes the upper second derivative of a 
function $f$ at $x_0$, as defined in \eqref{eq:uppersecondderivative}.
 
\begin{theorem}
\label{prop:inequality}
Let $\Om$ be a smooth convex body of a Riemannian manifold $(M^{n+1},g)$. Suppo\-se that the 
Ricci curvature of $M$ satisfies $\ric\geq n\delta$ on $\Om$. Then, the renormalized isoperimetric profile $Y_\Om=I_\Om^{(n+1)/n}$ verifies
\begin{equation}
\label{eq:inequality}
\overline{D^2Y_\Om}(V)\leq -(n+1)\,\delta\,Y_\Om(V)^{(1-n)/(1+n)}, \qquad V\in (0,\vol(\Om)).
\end{equation}

If equality holds for some $V_0\in (0,\vol(\Om))$ then the boundary $\Lambda=\overline{\ptl E\cap\Om}$ of any minimizer $E$ in $\Om$ of volume $V_0$ is a smooth, totally umbilical hypersurface such that
\[
\emph{Ric}(N,N)\equiv n\delta \ \text{ on } \Lambda\quad\text{ and }\quad \emph{II}(N,N)\equiv 0 \ \text{ on } \Lambda\cap\ptl\Om,
\]
where $N$ is the unit normal to $\Lambda$ which points into $E$, and \emph{II} is the second fundamental form of $\ptl\Om$ with respect to the inner normal.

Moreover, if $\Om$ coincides with the half-space $\mathbb{H}^{n+1}_\delta$ in the simply connected space form with constant sectional curvature $\delta$, then equality holds in \eqref{eq:inequality} for any $V\in (0,\vol(\Om))$.
\end{theorem}

\begin{proof}
Fix $V_0\in (0,\vol(\Om))$. By Proposition \ref{prop:regularity} there is a minimizer $E$ in $\Om$ of volume $V_0$. By the same result, the regular part $\Sg$ of $\Lambda=\overline{\ptl E\cap\Om}$ is a smooth, embedded hypersurface which meets $\ptl\Om$ orthogonally. The mean curvature of $\Sg$ with respect to the unit normal $N$ pointing into $E$ is a constant $H_0$. The boundary $\Sg\cap\ptl\Om$ could be empty, see Remark \ref{re:voidboundary}. In this case, we adopt the convention that the integrals over $\Sg\cap\ptl\Om$ are all equal to $0$. 

Consider a sequence of functions $\{\varphi_\eps\}_{\eps>0}$ as obtained in 
Lemma \ref{lem:phieps}. Fix $\eps>0$ and take a smooth vector field $X_\eps$ with compact 
support over $M$, such that $X_\eps(q)\in T_q (\ptl\Om)$ whenever $q\in\ptl\Om$ and 
$X_\eps=\varphi_\eps N$ in $\Sg$. The flow of diffeomorphisms 
$\{\phi_t\}_{t\in (-\gamma,\gamma)}$ of $X_\eps$ in $\overline{\Om}$ induces a variation $\{E_t=\phi_t(E)\}_t$ of $E$ through sets of finite perimeter contained in $\Om$. Call $\pp_\eps(t)=\pp(E_t,\Om)$ and $V_\eps(t)=\vol(E_t)$. By the first 
variation for perimeter and volume
\begin{align}
\label{eq:dp/dt}
\pp_\eps'(0)&=\int_\Sg\text{div}_\Sg\, X_\eps\,d\h_n=-\int_\Sg nH_0\,\varphi_\eps\,d\h_n, 
\\
\label{eq:dv/dt}
V_\eps'(0)&=\int_E\text{div}\,X_\eps\,d\h_{n+1}=-\int_\Sg\varphi_\eps\,d\h_n,
\end{align}
where $\text{div}_\Sg$ is the divergence relative to $\Sg$. As $V_\eps'(0)<0$, we can write $t$ as a function of the volume $V=V(t)$ for $V$ close to $V_0$; hence, we can define $\pp_\eps(V)=\pp_\eps[t(V)]$.

Now, consider the function $g_\eps(V)=\pp_\eps(V)^{(n+1)/n}$ defined on a neighborhood 
of $V_0$. By using the definition of isoperimetric 
profile and the fact that $E$ is a minimizer, it is clear that 
\[
Y_\Om(V)\leq g_\eps(V),\qquad Y_\Om(V_0)=g_\eps(V_0),
\]
from which we deduce
\begin{equation}
\label{eq:d2Y<=d2g}
\overline{D^2Y_\Om}\,(V_0)\leq\overline{D^2g_\eps}\,(V_0)=\left (\frac{n+1}{n}\right )\,\pp_\eps(V_0)^{1/n}\,
\left\{\frac{1}{n}\,\pp_\eps'(V_0)^2\,\pp_\eps(V_0)^{-1}+\pp_\eps''(V_0)\right\}.
\end{equation}

Now, we shall compute the derivatives $\pp'_\eps(V_0)$ and $\pp''_\eps(V_0)$. 
The first one is calculated by using \eqref{eq:dp/dt} and \eqref{eq:dv/dt}. We get 
\begin{equation}
\label{eq:dpeps/dv}
\pp_\eps'(V_0)=\pp_\eps'(0)\,V_\eps'(0)^{-1}=nH_0.
\end{equation}  

On the other hand, the calculation of $\pp_\eps''(V_0)$ requires second variation of perimeter and volume. By following the proof of \cite[Theorem 2.5]{zumbrun2} (in fact, the only change is that a new term involving the Ricci curvature appears), it is obtained  
\begin{align}
\label{eq:secondvariation}
\pp_\eps''(V_0)&=\left (\int_\Sg\varphi_\eps\,d\h_n\right )^{-2}
\\
\nonumber
&\times\left [\int_\Sg (|\nabla\varphi_\eps|^2-(\ric(N,N)+|\sigma|^2)\,
\varphi^2_\eps )\,d\h_n
-\int_{\Sg\cap\ptl\Om}\text{II}(N,N)\,\varphi_\eps^2\,d\h_{n-1}\right], 
\end{align}
where $|\sigma|^2$ is the squared sum of the principal curvatures of $\Sg$ with respect to $N$, and $\text{II}$ is the second fundamental form of $\ptl\Om$ with respect to the 
inner normal. 

Now, if we take $\limsup$ in the equality above when $\eps\to 0$ and 
we use Lemma \ref{lem:phieps} together with Fatou's Lemma, we have
\begin{align}
\label{eq:limiteps}
\limsup_{\eps\to 0}\,\pp_\eps''(V_0)&\leq-\pp(E,\Om)^{-2}\,
\biggl[\int_\Sg\big(\text{Ric}(N,N)+|\sg|^2\big)\,d\h_n
+\int_{\Sg\cap\ptl\Om}\text{II}(N,N)\,d\h_{n-1}\bigg]
\\
\nonumber
&\leq -n\,(\delta+H_0^2)\,\pp(E,\Om)^{-1},
\end{align}
where in the last inequality we have used the assumption on the Ricci curvature, the 
well-known inequality $|\sigma|^2\geq nH_0^2$, and the convexity of $\Om$.

Thus, if we pass to the limit in \eqref{eq:d2Y<=d2g} and we use \eqref{eq:dpeps/dv} together 
with \eqref{eq:limiteps}, we deduce
\begin{align}
\label{eq:twoinequalities}
\overline{D^2Y_\Om}(V_0)&\leq \biggl(\frac{n+1}{n}\biggr)\,\pp(E,\Om)^{1/n}\,
\biggl\{nH_0^2\,\pp(E,\Om)^{-1}
+\limsup_{\eps\to 0}\,\pp_\eps''(V_0)\biggr\}
\\
\nonumber
\notag
&\leq -(n+1)\,\delta\,\pp(E,\Om)^{(1-n)/n}=-(n+1)\,\delta\,Y_\Om(V_0)^{(1-n)/(1+n)},
\end{align}
and \eqref{eq:inequality} is proved. Moreover, if equality holds in \eqref{eq:twoinequalities} then we also have equality in \eqref{eq:limiteps}, and so $\Sg$ is totally umbilical, 
$\ric(N,N)\equiv n\delta$ on $\Sg$, and $\text{II}(N,N)\equiv 0$ on $\Sg\cap\ptl\Om$. Fur\-thermore, the singular set $\Sg_0=\Lambda-\Sg$ is empty by Proposition~\ref{prop:regularity} (iv) since $|\sg|^2$ is bounded. 
                                             
Finally, suppose that $\Om=\mathbb{H}^{n+1}_\delta$. By reflecting with respect to $\ptl\mathbb{H}^{n+1}_\delta$ we get that any minimizer in $\mathbb{H}^{n+1}_\delta$ is obtained by intersecting a geodesic ball $B$ centered at $\ptl\mathbb{H}^{n+1}_\delta$ with $\mathbb{H}^{n+1}_\delta$. As $\ptl B$ is a totally umbilical hypersurface and $\ptl\mathbb{H}^{n+1}_\delta$ is totally geodesic, we have equality in \eqref{eq:limiteps}. On the other hand, equality holds in the first inequality of \eqref{eq:twoinequalities} since $Y_{\mathbb{H}^{n+1}_\delta}$ equals the renormalized profile $\pp(V)^{(n+1)/n}$ given by parallel hypersurfaces to $\ptl B\cap\mathbb{H}^{n+1}_\delta$.
\end{proof}

\begin{remark}
\label{re:yOm}
By using the profile $h_\Om$ defined in \eqref{rmqopr} we easily see that Theorem \ref{prop:inequality} is also valid for the renormalized profile $y_\Om=h_\Om^{(n+1)/n}$. 
In particular,
\begin{equation}
\label{eq:inequality2}
\overline{D^2y_\Om}(\beta)\leq -(n+1)\,\delta\,y_\Om(\beta)^{(1-n)/(1+n)},\qquad\beta\in (0,1),
\end{equation}
with equality for all $\beta\in (0,1)$ when $\Om$ coincides with 
$\mathbb{H}^{n+1}_\delta$ $(\delta>0)$.
\end{remark}

\begin{remark}
Let $M$ be a closed Riemannian manifold with $\ric\geq n\delta$. Then, 
$M$ can be seen as a convex body $\Om$ with $\ptl\Om=\emptyset$, and 
the proof of \eqref{eq:inequality} remains valid with the only change 
that the terms involving $\Sg\cap\ptl\Om$ vanish. With a similar proof, 
V.~Bayle \cite[Theorem 1.1]{Ba} proved that \eqref{eq:inequality} holds 
for the renormalized profile $y_M=h_M^{(n+1)/n}$. Another type of 
differential inequality for the isoperimetric profile $I_M$ was previously 
established by F.~Morgan and D.~Johnson \cite[Proposition 3.3]{mj}.  
\end{remark}

The remainder of this section is devoted to deduce some immediate consequences of Theorem \ref{prop:inequality}. 

One of the easiest and most obvious applications of the differential inequality (\ref{eq:inequality}) is the following theorem. The proof only uses the fact that a lower semicontinuous function on an interval $f:I\to\rr$ such that $\overline{D^2f}\leq 0$ in the interior of $I$ must be concave (\cite{Bat}).

\begin{theorem}
\label{thconpc}
Let $\Omega$ be a smooth convex body of a Riemannian manifold $(M^{n+1},g)$ with non-negative 
Ricci curvature. Then, the renormalized profile of $\Om$ is concave. As consequence, the 
isoperimetric profile $I_\Om$ is concave and, therefore, increasing on 
$[0,\vol(\Omega)/2]$. 
\end{theorem}

\begin{remark}
The concavity of the profile of a smooth convex body 
$\Om\subset\rr^{n+1}$ was obtained by P.~Sternberg and 
K.~Zumbrun \cite[Theorem 2.8]{zumbrun2}. The observation 
that, in fact, the renormalized profile of 
$\Om\subset\rr^{n+1}$ is concave is due to E.~Kuwert \cite{K}.  
\end{remark}

Now, we generalize to the setting of convex bodies the regularity 
properties obtained for the isoperimetric profile of a closed 
Riemannian manifold (see \cite{BP}, \cite{mj} and \cite{Ba}). As 
an analytic outcome of Theorem~\ref{thconpc} we have that, under 
non-negative Ricci curvature, the isoperimetric profile of $\Om$ 
has the regularity properties of concave functions. In the 
following proposition we show that no assumption on the Ricci 
curvature is needed.  

\begin{proposition}
\label{syprdn}
Let $\Om$ be a smooth convex body of a Riemannian manifold $(M^{n+1},g)$. Then 
the renormalized isoperimetric profile $Y_\Om$ has left and right derivatives satisfying 
$$(Y_\Om)'_l(V)\geq (Y_\Om)'_r(V),\quad V\in (0,\vol(\Om)).$$
As a consequence, $(Y_\Om)'_l (\vol(\Omega)/2)$ is non-negative. Moreover, $Y_{\Omega}$ is differentiable on the interval $(0,\vol(\Omega))$ except on an at most countable set. 

As to the isoperimetric profile $I_\Om$, it has left and right derivatives for every $V$ in $\big(0,\vol(\Omega)\big)$, such that
$$(I_\Om)_l'(V)\geq nH_E\geq(I_\Om)_r'(V),$$
where $H_E$ is the inward mean curvature associated to a minimizer $E$ in $\Om$ of volume $V$.
As a consequence, $(I_\Om)_l' (\vol(\Omega)/2)$ is non-negative. Furthermore, $I_\Om$ is differentiable on $\big(0,\vol(\Omega)\big)$ except on an at most countable set. 
\end{proposition}

\begin{proof}
The regularity properties and the inequality between the side derivatives come from the differential inequality (\ref{eq:inequality}), which implies that locally around $V_0\in (0,\vol(\Om))$ the renormalized profile $Y_\Om$ is concave, up to the addition of a constant times $(V-V_0)^2$. Now fix $V_0\in (0,\vol(\Om))$ and take a minimizer $E$ in $\Om$ of volume $V_0$. Let $\pp(V)$ be the relative profile associated to an almost parallel variation of $E$ as constructed in the proof of Theorem \ref{prop:inequality}. It is obvious that $I_\Om(V)\leq\pp(V)$ for $V$ close to $V_0$, and $I_\Om(V_0)=\pp(V_0)$. As $\pp'(V_0)=nH_E$ (see \eqref{eq:dpeps/dv}) we deduce that 
$(I_\Om)'_l(V_0)\geq nH_E\geq (I_\Om)'_r(V_0)$. 
\end{proof}

\begin{remark}
\label{ancontmc}
The asymptotic behaviour of the profile given in Proposition \ref{asbpnz} allows us to deduce the following consequences from Proposition \ref{syprdn}
\begin{itemize}
\item[(i)] $I_\Om$ is continuous on $[0,\vol(\Om)]$.
\item[(ii)] $\lim_{V\to 0}\, (I_\Om)'_r(V)=+\infty$.
\item[(iii)] The inward mean curvature associated to a minimizer in $\Om$ explodes when the enclosed volume tends to zero.  
\end{itemize} 
\end{remark}

We finish this section by showing some topological restrictions related 
to the connectivity of minimizers inside a convex body. We derive them 
by using a well-known argument (see \cite{zumbrun2} and \cite{mj}), which 
relies on the second variation formula of perimeter \eqref{eq:secondvariation}.

\begin{proposition}
\label{syprdn2}
Let $\Om$ be a smooth convex body of a Riemannian manifold $(M^{n+1},g)$ such that 
$\ric\geq n\delta$ on $\Om$. Denote by \emph{II} the second fundamental 
form of $\ptl\Om$ with respect to the inner normal. Let $E$ be an isoperimetric region 
in $\Om$, $\Sg$ the regular part of $\overline{\ptl E\cap\Om}$, and $N$ the normal 
to $\Sg$ pointing into $E$. Then
\begin{itemize}
\item[(i)] If $\delta>0$, then $\Sg$ is connected. 
\item[(ii)] If $\delta=0$ and $\Sg$ consists of more than one component, then $\Sg$ is totally geodesic and we have $\ric(N,N)\equiv 0$ in $\Sg$ and $\emph{II}(N,N)\equiv 0$ in 
$\Sg\cap\ptl\Om$. As consequence, if $\Sg$ is non-connected, and $\Om$ is strictly convex 
in the sense that $\emph{II}>0$, then $\Sg\cap\ptl\Om=\emptyset$.  
\item[(iii)] If $\delta\leq 0$, then there exists $V_1\in (0,\vol(\Om))$ 
such that $\Sg$ is connected if $\vol(E)\leq V_1$.   
\end{itemize}
\end{proposition}
  
\begin{proof}
Call $V_0=\vol(E)$ and denote by $H_0$ the mean curvature of $\Sg$ with respect to $N$. Let $\Sg'$ be a component of $\Sg$ and $\{\varphi_\eps\}\subset C^\infty_0(\Sg')$ a sequence as in Lemma \ref{lem:phieps}. By following the proof of Theorem \ref{prop:inequality} we consider almost parallel variations of $E$ and the associated perimeter functions $\pp_\eps(V)$. Call $\alpha(V_0)=\limsup_{\eps\to 0}\pp''_\eps(V_0)$. From \eqref{eq:limiteps} we know that
\begin{flushleft}
(*) \ \qquad\qquad\qquad\qquad\qquad\qquad $\alpha(V_0)\leq -n\,(\delta+H_0^2)\,\pp(E,\Om)^{-1}$,
\end{flushleft}
due to the hypothesis on the Ricci curvature, the convexity of $\Om$ and the inequality $|\sg|^2\geq nH_0^2$.

We assert that $\alpha(V_0)<0$ implies that $\Sg$ is connected. Otherwise, we would use almost parallel variations with $\eps\approx 0$ to expand one component $\Sg_1$ and shrink another one $\Sg_2$ so that the resulting variation preserves volume while reducing perimeter, see 
\cite[Theorem 2.6]{zumbrun2} for details; this would give us a contradiction with the minimality of $E$.

Now we distinguish two cases. If $\delta\geq 0$, then $\alpha(V_0)\leq 0$ and an easy 
discussion of equality cases in (*) proves the claim. If $\delta\leq 0$, then the 
explosion of the mean curvature~for small volumes (Remark \ref{ancontmc} (iii)) yields 
the existence of $V_1$ such that $\alpha(V)<0$ for $V\in [0,V_1]$.    
\end{proof}

\begin{remark}
Topological restrictions on isoperimetric hypersurfaces inside a Euclidean 
convex body were obtained by A.~Ros and E.~Vergasta \cite{rosvergasta} and by 
P.~Sternberg and K.~Zumbrun \cite{zumbrun2}. On the one hand, statement (ii) 
in the proposition above is proved in \cite[Theorem 2.6]{zumbrun2} for a convex 
body $\Om\subset\rr^{n+1}$. Furthermore, it is shown that strict convexity of $\Om$ 
implies that $\Sg$ is connected. We must point out that, in general, this cannot 
be achieved when $\Om$ is not a Euclidean domain since $\Sg\cap\ptl\Om$ could be 
empty, see Remark \ref{re:voidboundary}. On the other hand, in 
\cite[Theorem 5]{rosvergasta} some conditions on the genus $g$ and the 
number $r$ of boundary components of $\Sg$ are obtained when $\Om\subset\rr^3$. 
In precise terms, they proved that the only possible values for $g$ and $r$ are
\begin{itemize}
\item[(i)] $g=0$ and $r=1,2$ or $3$;
\item[(ii)] $g=2$ or $3$ and $r=1$.
\end{itemize} 
It has been recently conjectured that an isoperimetric hypersurface inside a 
strictly convex body of $\rr^3$ must be homeomorphic to a disk (\cite{ros}).
\end{remark}

Let $\Om$ be a smooth convex body of a Riemannian manifold and let $n\delta$ be a lower bound on the Ricci curvature of $\Om$. By Proposition \ref{syprdn2} we have that a minimizer $E$ in $\Om$ is connected when $\delta>0$, or when $\delta\leq 0$ and $\vol(E)$ is small enough. At first, the second variation of perimeter is not sufficient, in the case $\delta\leq 0$, to discard a minimizer with finitely many components bounded by totally geodesic hypersurfaces. However, by using that the profile is concave when $\delta=0$ we can prove

\begin{proposition}
\label{prop:syprdn}
Let $\Om$ be a smooth convex body of a Riemannian manifold with non-negative Ricci curvature. Then, isoperimetric regions in $\Om$ are connected.
\end{proposition}

\begin{proof}
Suppose that $E$ is a minimizer of volume $V_0\in (0,\vol(\Om))$ and that 
$E_1$ is a connected component of $E$ with volume $V_1<V_0$. By the definition of isoperimetric profile and the fact that the set of singularities in $\overline{\ptl E\cap\Om}$ does not contribute to perimeter, we get 
\[
I_\Om(V_0)=\pp(E,\Om)=\pp(E_1,\Om)+\pp(E-E_1,\Om)\geq I_\Om(V_1)+I_\Om(V_0-V_1).
\]
On the other hand, the concavity of $Y_{\Omega}$ (Theorem \ref{thconpc}) gives us
\[
Y_\Om(V_0)\leq Y_\Om(V_1)+Y_\Om(V_0-V_1),
\]
and so, as $I_{\Omega}(V_1)$ and $I_\Om(V_0-V_1)$ are positive, and since the 
function $x\longmapsto x^{\frac{n}{n+1}}$ is strictly concave, we deduce 
\[
I_\Om(V_0)<I_\Om(V_1)+I_\Om(V_0-V_1),
\]
which leads us to a contradiction. This proves that $V_1=V_0$, and $E$ is therefore connected.
\end{proof}

\section{Comparison theorems}
\label{compmod}
\setcounter{equation}{0}

In this section, we shall integrate the differential inequality 
\eqref{eq:inequality} in order to prove comparison theorems for the 
isoperimetric profile of a smooth convex body $\Om$ in a Riemannian 
manifold $(M^{n+1},g)$. The underlying philosophy of these results 
consists in using the analytic Theorem \ref{th:analcomp} to compare 
a profile $f$, which can be $Y_\Om$ or the function $y_\Om$ defined 
in Remark \ref{re:yOm}, with a solution $g$ of the differential 
equation associated to \eqref{eq:inequality} having the same 
initial conditions or the same boundary values as $f$. In the first 
case we shall obtain an upper bound for the profile $I_\Om$, while in 
the second one, we shall deduce a lower bound for $h_\Om$ that can be 
interpreted as a L\'evy-Gromov type inequality. We must remark that 
both comparisons are quite different although they arise from the same 
differential inequality. A detailed analysis of equality cases will 
allow us to deduce global geometric consequences on $\Om$.  

Through this section we also illustrate how to use our isoperimetric 
inequalities to deduce other geometric an analytic comparisons. In 
this way, we give alternative proofs of the inequalities 
in~Theorem~\ref{th:bishop}, and we characterize the half-spheres as 
the only convex domains for which equality in 
Theorem~\ref{th:bishop} (i) holds. Finally, we prove a comparison result 
for the first non-zero Neumann eigenvalue of the Laplace operator on 
$\Om$ that can be seen as a generalization of the Obata--Lichnerowicz 
theorem~\cite[Theorem 9, p. 82]{chavel3}.

\subsection{Upper bounds on the isoperimetric profile}

\begin{theorem}
\label{thcomp}
Let $\Omega$ be a smooth convex body with non-empty boundary of a complete Riemannian 
manifold $(M^{n+1},g)$. Suppose that the Ricci curvature of $M$ satisfies $\ric\geq n\delta$ on $\Om$. Then 
\begin{equation}
\label{incoprchp}
I_{\Omega}(V)\leq I_{\mathbb{H}_{\delta}^{n+1}}(V), \qquad V\in [0,\vol(\Om)],
\end{equation}
where $\mathbb{H}^{n+1}_\delta$ is a half-space in the $(n+1)$-dimensional simply connected space form with constant sectional curvature $\delta$.

If equality holds in \eqref{incoprchp} for some $V_0\in (0,\vol(\Omega)]$, then  $I_{\Omega}= I_{\mathbb{H}_{\delta}^{n+1}}$ on $[0,V_0]$, and the boundary $\overline{\ptl E\cap\Om}$ of any minimizer $E$ in $\Om$ of volume $V\in (0,V_0)$ is a smooth, totally umbilical hypersurface. Moreover, if $V_0=\vol(\Om)$ \emph{(}which implies $\delta>0$\emph{)} then $\Om$ is 
isometric to~$\mathbb{H}^{n+1}_\delta$.
\end{theorem}

\begin{proof}
The comparison arises from the fact that a continuous solution of the differential 
inequality \eqref{eq:inequality} is bounded from above by a solution of the differential 
equation 
\begin{equation}
\label{eq:f2}
f''=-(n+1)\,\delta\,f^{(1-n)/(1+n)}
\end{equation}
with the same initial conditions (Theorem \ref{th:analcomp}). Then, by using  
that the renormalized profile of $\mathbb{H}^{n+1}_\delta$ satisfies \eqref{eq:f2} (see Theorem \ref{prop:inequality}) 
and taking into account the asymptotic behaviour of the renormalized profile $Y_\Om$ at the origin (Proposition \ref{asbpnz}), we obtain
\begin{equation}
\label{eq:hola}
Y_\Om(V)\leq Y_{\mathbb{H}^{n+1}_\delta}(V), \qquad V\in [0,\min\{\vol(\Om),\vol(\mathbb{H}_\delta^{n+1})\}]. 
\end{equation}
From the inequality above we get \eqref{incoprchp} once we show that 
$\vol(\Om)\leq\vol(\hh_\delta^{n+1})$. This volume comparison is trivial
if $\delta\leq 0$ while in the case $\delta>0$, the opposite inequality
would allow us to deduce from \eqref{eq:hola} that 
$Y_\Om(\vol(\hh_\delta^{n+1}))\leq 0$, which is a contradiction since the
profile is positive in $(0,\vol(\Om))$. 
 
Finally, if both profiles coincide at $V_0\in (0,\vol(\Om)]$ then they
must coincide in $[0,V_0]$ by Theorem \ref{th:analcomp}. The umbilicality 
of a minimizer of volume $V<V_0$ follows from the discussion, given in 
Theorem \ref{prop:inequality}, of equality cases in \eqref{eq:inequality}. 
If $V_0=\vol(\Om)$ then $\vol(\Om)=\vol(\mathbb{H}^{n+1}_\delta)$ and $\Om$ 
is isometric to $\mathbb{H}^{n+1}_\delta$ by Theorem \ref{th:bishop} (ii). 
\end{proof}

\begin{remark}
\label{re:altproofca}
Note that we have given another proof of the volume comparison 
$\vol(\Om)\leq\vol(\hh_\delta^{n+1})$ of Theorem~\ref{th:bishop} (ii) by 
using the isoperimetric inequality \eqref{eq:hola}.
\end{remark}

\begin{remark}
When $n=1$ the differential inequality \eqref{eq:inequality} turns out 
to be linear and Theorem \ref{incoprchp} follows since the function 
$E(V)=Y_\Om(V)-Y_{\mathbb{H}^2_\delta}(V)$ is concave on $[0,\vol(\Om)]$ 
and the tangent line at the origin coincides with the $x$-axis. After 
an explicit calculation of $Y_{\hh^2_\delta}$, inequality \eqref{incoprchp} 
reads
\[
I^2_\Om(V)\leq V\,(2\pi-\delta V),\qquad V\in [0,\vol(\Om)].
\]
\end{remark}

\begin{remark}
For a closed Riemannian manifold $M$ with 
$\ric\geq n\delta$ the integration of the differential inequality 
\eqref{eq:inequality} would give us the comparison
\begin{equation}
\label{eq:noboundary}
I_M\leq I_{\mathbb{M}^{n+1}_\delta},\qquad V\in [0,\vol(M)],
\end{equation}
where $\mathbb{M}^{n+1}_\delta$ stands for the $(n+1)$-dimensional 
simply connected space form with constant sectional curvature 
$\delta$. This result was previously proved by F.~Morgan and D.~Johnson 
\cite[Theorem 3.4]{mj}. 
\end{remark}

\begin{remark}
\label{re:mj}
Inequality \eqref{incoprchp} is also valid for a smooth, unbounded, 
convex domain $\Om$ with non-empty boundary and $\ric\geq n\delta$. 
This can be proved by showing, as was done in \cite[Theorem 3.5]{mj} 
for closed Riemannian manifolds, that the perimeter in $\Om$ of a 
``half-ball'' $B=\Om\cap B(p,r)$ centered at a point $p\in\ptl\Om$ 
is less than or equal to the area of the geodesic half-ball $\widetilde{B}$ 
in $\hh^{n+1}_\delta$ of the same volume, with equality only if $B$ 
is isometric to $\widetilde{B}$ and $\ptl\Om$ is geodesic at $p$. 
The arguments in the proof by F.~Morgan and D.~Johnson rely on  
comparison theorems involving the volume of metric balls 
(\cite[Theorem 3.9]{chavel}) and the area of metric spheres 
(\cite[Proposition 3.3]{chavel}). These theorems do not use the 
compactness of the ambient manifold and can be generalized to our setting 
by following the scheme in the proof of Theorem~\ref{th:bishop}.

This alternative proof of \eqref{incoprchp} also allows us to deduce 
geometric consequences on $\Om$ when we have equality in \eqref{incoprchp}. 
We summarize them in the next result 

\begin{theorem}
\label{th:unbounded}
Let $\Om$ be a smooth convex domain with $\ric\geq n\delta$ in a complete Riemannian
manifold $(M^{n+1},g)$. Then 
\begin{itemize}
\item[(i)] If $\Om$ has non-empty boundary then \eqref{incoprchp} holds, and equality 
for some $V_0\in (0,\vol(\Om)]$ implies that $\ptl\Om$ is totally geodesic in $M$ and 
$\Om$ has constant sectional curvature $\delta$ in a neighborhood of $\ptl\Om$. 
\item[(ii)] If $\ptl\Om=\emptyset$ then \eqref{eq:noboundary} holds, and equality for 
some $V_0\in (0,\vol(\Om)]$ implies that $M$ is isometric to a quotient of the 
simply connected space form $\mathbb{M}^{n+1}_\delta$ with constant sectional 
curvature $\delta$.
\end{itemize} 
\end{theorem}    
\end{remark}

\begin{remark}
In general, we cannot improve statement (i) in the theorem above to the stronger 
conclusion that equality in \eqref{incoprchp} for some $V_0$ implies that $\Om$ 
has constant sectional curvature $\delta$. For example, denote by $\Om$ the domain 
obtained from attaching the half-sphere of $\mathbb{S}^2$ centered at the north pole to 
the compact cylinder $\mathbb{S}^1\times [-1,0]$ through the circle 
$\mathbb{S}^1\times\{0\}$. It is clear that $I_\Om=I_{\hh^2_0}$ for small 
values; however, $\Om$ is not a flat domain.
\end{remark}

\subsection{A L\'evy-Gromov type inequality for convex bodies}

Let $(M^{n+1},g)$ be a closed Riemannian manifold with $\ric\geq n\delta>0$. Denote 
by $h_M$ the profile of $M$ as defined in \eqref{rmqopr}. 
L\'evy-Gromov inequality \cite{GR2} states that
\begin{equation}
\label{eq:lg}
h_M(\beta)\geq h_{\mathbb{M}^{n+1}_\delta}(\beta),\qquad\beta\in [0,1],
\end{equation}
where $\mathbb{M}^{n+1}_\delta$ is an $(n+1)$-dimensional sphere of radius 
$1/\sqrt{\delta}$. Moreover, if equality holds in \eqref{eq:lg} for some 
$\beta\in (0,1)$, then $M$ is isometric to $\mathbb{M}^{n+1}_\delta$.  

Inequality \eqref{eq:lg} can be obtained by integrating a differential 
inequality similar to (\ref{eq:inequality}), see \cite{Bat}. With a 
similar technique, we generalize \eqref{eq:lg} to the setting of convex bodies. 

\begin{theorem}
\label{generallevygromov}
Let $\Om$ be a smooth convex body of a Riemannian manifold $(M^{n+1},g)$. Suppose 
that the Ricci curvature of $M$ over $\Om$ satisfies $\ric\geq n\delta>0$. Then,
\begin{equation}
\label{eq:ourlg}
h_{\Omega}(\beta)\geq h_{\mathbb{H}^{n+1}_\delta}(\beta),\qquad\beta\in [0,1], 
\end{equation}
where $\mathbb{H}^{n+1}_\delta$ is an $(n+1)$-dimensional half-sphere of radius $1/\sqrt{\delta}$.  

Moreover, if $\Om$ has non-empty boundary then equality holds in 
\eqref{eq:ourlg} for some $\beta_0\in (0,1)$ if and only if $\Om$ 
is isometric to $\mathbb{H}^{n+1}_\delta$. 
\end{theorem}

\begin{proof}
The inequality follows from the fact, given in Theorem 
\ref{th:analcomp} (i), that any function satisfying the 
differential inequality \eqref{eq:inequality2} is bounded 
from below by an exact solution of the corresponding 
differential equation with the same boundary values. 
Furthermore, if we have equality for some $\beta_0\in (0,1)$ 
then $h_\Om=h_{\mathbb{H}^{n+1}_\delta}$ 
on $[0,1]$, and by the asymptotic behaviour of $h_\Om$ at 
the origin (Proposition \ref{asbpnz}) we deduce that 
$\vol(\Om)=\vol(\mathbb{H}^{n+1}_\delta)$. From statement (ii) 
in Theorem \ref{th:bishop} we conclude that~$\Om$ is isometric 
to $\mathbb{H}^{n+1}_\delta$. 
\end{proof}

The preceding result can be given in the following alternative form
 
\begin{corollary}
\label{cgenerallevygromov}
Let $\Om$ be a smooth convex body of a Riemannian manifold $(M^{n+1},g)$. Suppose
that the Ricci curvature of $M$ over $\Om$ satisfies $\ric\geq n\delta>0$. 
Then, for any Borel set $E\subseteq\Om$, we have
\[
\frac{\mathcal{P}(E,\Omega)}{\vol(\Omega)}\geq \frac{\mathcal{P}(E^*,\mathbb{H}^{n+1}_{\delta})}{\vol(\mathbb{H}^{n+1}_{\delta})},
\]
where $E^*\subseteq\mathbb{H}^{n+1}_\delta$ is a geodesic half-ball centered at 
$\ptl\mathbb{H}^{n+1}_{\delta}$ such that 
$$\frac{\vol(E)}{\vol(\Omega)}=\frac{\vol(E^*)}{\vol(\mathbb{H}^{n+1}_{\delta})}.$$

Moreover, if $\Om$ has non-empty boundary and equality holds for some set 
$E\subseteq\Om$ with $\vol(E)\in (0,\vol(\Om))$, then $\Om$ is isometric 
to an $(n+1)$-dimensional half-sphere of radius~$1/\sqrt{\delta}$. 
\end{corollary}

\begin{remark}
\label{rem:cheegerconstant}
Let $h_C(\Om)$ be the Cheeger isoperimetric constant of a smooth convex body 
$\Omega$ of a Riemannian manifold $(M^{n+1},g)$, defined by
\[
h_C(\Om)=\inf\bigg\{\frac{\pp(E,\Om)}
{\min\,\{\vol(E),\vol(\Om\setminus E)\}}:\vol(E)\in (0,\vol(\Om))\bigg\}.
\]
Note that
\[
h_C(\Om)=\inf\bigg\{\frac{h_\Om(\beta)}
{\min\,\{\beta,1-\beta\}}:\beta\in (0,1)\bigg\},
\]
and so, if the Ricci curvature of $M$ is non-negative on $\Omega$, we deduce 
by the concavity of the profile (Theorem~\ref{thconpc})
\[
h_C(\Om)=2\,h_{\Om}(1/2),
\]
which yields $h_C(\Om)\geq h_C(\hh_\delta^{n+1})$ when $\ric\geq n\delta>0$ in $\Om$ 
by \eqref{eq:ourlg}.

Now, by reproducing the arguments in \cite{Ba} (see also \cite{Bat}), we can 
refine Theorem \ref{generallevygromov}, so as to get, under the same  
assumption on the Ricci curvature,
\begin{eqnarray}
\label{levygromovrefined}
h_{\Omega}(\beta)\geq \bigg[\frac{h_C(\Omega)}{h_C(\mathbb{H}^{n+1}_{\delta})}\bigg]^{\frac{1}{n+1}} h_{\mathbb{H}^{n+1}_{\delta}}(\beta),\qquad\beta\in[0,1].
\end{eqnarray}  
Moreover, if there is $\beta_0\in (0,1)$ such that (\ref{levygromovrefined}) 
is an equality, then $\Omega$ is isometric to $\mathbb{H}^{n+1}_{\delta}$.
\end{remark}

\subsection{Some consequences of Theorem~\ref{generallevygromov}}
\label{eq:riemanniansetting}

We first show how to use Theorem~\ref{generallevygromov} to give a characterization 
of equality cases in Theorem~\ref{th:bishop} (i). We need a previous result 
(see \cite{G1} for closed Riemannian manifolds), linking the diameter of a domain 
$\Om$ and the profile $h_\Om$.

\begin{lemma}
\label{lem:intineq}
The diameter of a smooth domain $\Om$ of a complete Riemannian manifold $(M^{n+1},g)$ 
satisfies
\[
\emph{diam}(\Om)\leq\int_0^1\frac{d\beta}{h_\Om(\beta)},
\]
with equality when $\Om$ coincides with an $(n+1)$-dimensional half-sphere.
\end{lemma}
\begin{proof}

Suppose that $\vol(\Om)<\infty$ (in other case $h_\Om\equiv 0$). If $\Om$ 
is unbounded then choose any point 
$p_0\in\Om$. If $\Om$ is bounded, fix a point $p_0\in\overline{\Om}$ such 
that $\text{dist}(p_0,p_1)=\text{diam}(\Om)$ for some $p_1\in\overline{\Om}$. 
Denote by $S_t$ and $B_t$ the metric sphere and the metric 
open ball in $M$ centered at $p_0$ with radius $t>0$. By the coarea 
formula~\cite[Corollary I.3.1]{C2}, the volume of a set $E\subseteq M$ can be
computed as
\[
\vol(E)=\int_0^{+\infty}\h_n(E\cap S_t)\,dt,
\]
and so the function $\beta(r)=\vol(\Om\cap B_r)/\vol(\Om)$ is 
absolutely continuous on $[0,\text{diam}(\Om)]$ and satisfies
\begin{equation}
\label{eq:torillo}
\beta'(r)=\frac{\h_n(\Om\cap S_r)}{\vol(\Om)}\geq
\frac{\pp(\Om\cap B_r,\Om)}{\vol(\Om)}\geq
h_\Om(\beta(r)),
\end{equation}
for almost all $r\in [0,\text{diam}(\Om)]$, with equality when 
$\Om$ coincides with a half-sphere. The proof finishes by 
integrating in \eqref{eq:torillo}. 
\end{proof}

\begin{remark}
The asymptotic behaviour of $h_\Om$ at the origin (Proposition~\ref{asbpnz}) 
ensures that the upper bound on the diameter given in the lemma above 
is finite when $\Om$ is bounded.
\end{remark}

As a consequence of Lemma~\ref{lem:intineq} and 
Theorem~\ref{generallevygromov} we can prove for convex 
bodies the analogous of the well-known Topogonov--Cheng 
theorem \cite[Theorem 3.11]{chavel} for closed Riemannian 
manifolds. Note that the following result is not a direct 
consequence of the aforementioned one for closed manifolds 
since we are assuming that $\ric\geq n\delta>0$ 
only in $\Om$.
\begin{theorem}
\label{th:cheng}
Let $\Om$ be a smooth convex body with non-empty boundary of a Riemannian manifold
$(M^{n+1},g)$. If the Ricci curvature of $M$ satisfies $\ric\geq n\delta>0$ on 
$\Om$, then
\[
\emph{diam}(\Om)\leq\frac{\pi}{\sqrt{\delta}},
\]
and equality holds if and only if $\Om$ is 
isometric to a half-sphere of radius $1/\sqrt{\delta}$.
\end{theorem}

\begin{remark}
By following the arguments in \cite[Theorems 3.2 and 3.3]{Ba} we could say 
that, for a smooth convex body $\Om$ with non-empty boundary and 
$\ric\geq n\delta>0$, having a diameter close to $\pi/\sqrt{\delta}$ 
(resp. a volume close to $\vol(\hh_\delta^{n+1})$) is equivalent to the fact 
that $h_{\Omega}-h_{\mathbb{H}_{\delta}^{n+1}}$ is uniformly close to $0$ 
on $[0,1]$ (resp. $h_\Om/h_{\hh_\delta^{n+1}}$ is uniformly close to $1$ on 
$(0,1)$). This means that almost maximality of the diameter or almost 
maximality of the volume both entail, in certain sense, almost minimality of 
the profile.   
\end{remark}

We finish this section with an eigenvalues comparison theorem. The application 
of an isoperimetric inequality to obtain eigenvalues estimates
(see \cite[Theorem 2, p. 87]{chavel3}) was first given by G. B. Faber and 
E. Krahn for smooth Euclidean domains with compact closure. In \cite{bemeyer} and 
\cite{BBG} it is shown how the ideas of G. B. Faber and E. Krahn, together 
with L\'evy--Gromov inequality \eqref{eq:lg}, lead to sharp estimates for the 
first eigenvalue of the Laplace operator with Dirichlet boundary condition on a 
smooth, bounded domain of a complete Riemannian manifold $(M^{n+1},g)$ with 
$\ric\geq n\delta>0$. Other estimates for Dirichlet eigenvalues obtained in a 
similar way can be found in \cite{G1} and \cite{Ba}. 

In the setting of a smooth convex domain $\Om$ with $\ptl\Om\neq\emptyset$, 
the fact that isoperimetric hypersurfaces in the model $\hh_\delta^{n+1}$ intersect 
the boundary orthogonally, seems to indicate that the Neumann boundary condition 
on $\ptl\Om$ is more appropriated if we want to derive an eigenvalues comparison 
from inequality \eqref{eq:ourlg}. In fact, we can prove 

\begin{theorem}
\label{theocompeigenvalue}
Let $\Om$ be a smooth convex body with non-empty boundary of a complete 
Riemannian manifold $(M^{n+1},g)$. If the Ricci curvature of $M$ 
satisfies $\ric\geq n\delta>0$ on $\Om$, then 
\begin{eqnarray}
\label{complambuneum}
\lambda_1^N(\Om)\geq\lambda_1^N(\hh^{n+1}_{\delta})=(n+1)\,\delta,
\end{eqnarray}
where the notation $\lambda_1^N(\Omega)$ stands for the lowest non-zero 
eigenvalue of the Laplace operator on $\Om$ with Neumann boundary condition 
on $\partial\Omega$. Moreover, if \eqref{complambuneum} is an equality, then 
$\Om$ is isometric to a half-sphere $\mathbb{H}^{n+1}_{\delta}$ of radius 
$1/\sqrt{\delta}$.
\end{theorem}

\begin{proof}
We give a brief desciption of the proof, which follows the symmetrization 
argument in \cite[Th\'eor\`eme 5]{bemeyer}. For any  
non-trivial function $u\in C^\infty(\Om)$, denote by $R_\Om(u)$ the Rayleigh 
quotient of $u$, given by
\[
R_\Om(u)=\bigg ( \int_\Om |\nabla u|^2\,d\h_{n+1}\bigg)\,\bigg (\int_\Om u^2\,d\h_{n+1} 
\bigg )^{-1}.
\]
Due to the variational characterization of Neumann eigenvalues there exists a 
smooth, mean zero function $u$ on $\overline\Om$ such that 
$R_\Om(u)=\lambda_1^N(\Om)$ and $\ptl u/\ptl\nu=0$ on $\ptl\Om$, where $\nu$ 
is the inward normal vector to $\ptl\Om$. Suppose that $u$ has finitely many 
non-degenerate critical points (condition (ND)). The symmetrization technique 
allows us to construct, by using a suitable family of concentric half-balls in 
$\hh_\delta^{n+1}$ centered at a fix boundary point, a function $u^*$ defined on $\hh_\delta^{n+1}$ such that
\begin{itemize}
\item[(i)] $u^*$ is a non-trivial Sobolev function on $\hh_\delta^{n+1}$.
\item[(ii)] $u^*$ has mean zero over $\hh_\delta^{n+1}$.
\item[(iii)] $R_\Om(u)\geq R_{\hh_\delta^{n+1}}(u^*)$ with equality if and only if 
$\Om$ is isometric to $\hh_\delta^{n+1}$ (here is the point where Theorem 
\eqref{generallevygromov} 
is used).
\end{itemize} 
By using statement (iii) and the variational characterization of Neumann 
eigenvalues, the proof of the theorem follows. 

If $u$ does not satisfy condition (ND), then we get \eqref{complambuneum} 
by approximation since $\lambda_1^N(\Om)$ is the limit of a sequence 
$\{R_\Om(u_n)\}_{n\in\nn}$, where each $u_n$ has mean zero and satisfies 
condition (ND). In this situation, the discussion of the equality case 
is not so obvious; we appeal to \cite[p. 520]{bemeyer}.
\end{proof}  

\begin{remark}
By using inequality \eqref{levygromovrefined} instead of \eqref{eq:ourlg} 
in the proof of Theorem~\ref{theocompeigenvalue}, we obtain
\[
\label{complambuneumint}
\lambda_1^N(\Omega)\geq\bigg[\frac{h_C(\Omega)}
{h_C(\mathbb{H}^{n+1}_{\delta})}\bigg]^{\frac{2}{n+1}} 
\lambda_1^N(\mathbb{H}^{n+1}_{\delta}),
\]
with equality if and only if $\Om$ is isometric to $\hh_\delta^{n+1}$.
\end{remark}

\section{Appendix: an alternative proof of inequality 
\eqref{incoprchp} in the euclidean case}
\label{sec:alternative}
\setcounter{equation}{0}

Here we give a geometric proof of the fact that the isoperimetric profile of a convex body $\Om\subset\rr^{n+1}$ is bounded from above by the profile of the half-space $\mathbb{H}^{n+1}=\{x_{n+1}>0\}$. The proof relies on the fact that the local convexity of a domain $\Om$ around a boundary point implies $I_\Om\leq I_{\mathbb{H}^{n+1}}$ for small volumes.

\begin{proposition}
\label{prop:alternative}
Let $\Om$ be a smooth domain in $\rr^{n+1}$. If $\Om$ has a local supporting hyperplane at a point $x\in\ptl\Om$, then there exists $V_0>0$ such that $I_\Om(V)\leq I_{\mathbb{H}^{n+1}}(V)$, whenever $V\in [0,V_0]$. 
\end{proposition}

\begin{proof}
We follow the proof in \cite[Proposition 3.6]{ritros}. Denote by $\pp(r)$ and $V(r)$ 
respectively the perimeter in $\Om$ and the volume of the ball $B_r$ of radius $r>0$ 
centered at $x$ intersected with $\Om$. Let $\wt{V}(r)$ be the volume of the cone 
subtended by $\ptl B_{r}\cap\Om$ and vertex at $x$. We have the relation
\[
\pp(r)=(n+1)\,\frac{\wt{V}(r)}{r}.
\]

On the one hand, since $\Om$ is locally convex around $x$, we have $V(r)\ge\wt{V}(r)$ for $r$
small, so that
\[
\pp(r)=(n+1)\,\frac{\wt{V}(r)}{r}\le (n+1)\,\frac{V(r)}{r}.
\]
On the other hand, if $\pp_{e}(r)$ and $V_{e}(r)$ respectively are the 
area and the volume of a half-ball in $\mathbb{H}^{n+1}$ of radius $r>0$, we have
\[
\frac{\pp_{e}(r)}{V_{e}(r)}=\frac{n+1}{r},
\]
and so
\[
\frac{\pp(r)}{V(r)}\le\frac{\pp_{e}(r)}{V_{e}(r)}.
\]
Since $V(r)\le V_{e}(r)$ due to the local convexity of $\Om$ around $x$, we finally get
\[
\frac{\pp(r)}{V(r)^{n/(n+1)}}=\frac{\pp(r)}{V(r)}\,V(r)^{1/(n+1)}\le
\frac{\pp_{e}(r)}{V_{e}(r)}\,V_{e}(r)^{1/(n+1)}=
\frac{\pp_{e}(r)}{V_{e}(r)^{n/(n+1)}}=d_{n},
\]
where $d_{n}$ is the constant that appears in the expression of the 
isoperimetric profile of the half-space $I_{\mathbb{H}^{n+1}}(V)= d_nV^{n/(n+1)}$.

Hence, for small $r$, we obtain the relation $\pp(r)\le I_{\mathbb{H}^{n+1}}(V(r))$, which 
proves the claim.
\end{proof}

\noindent
\emph{Proof of inequality \eqref{incoprchp}:} Let $\Om$ be a smooth convex 
body in $\rr^{n+1}$. As the renormalized profile $Y_{\mathbb{H}^{n+1}}$ is 
linear as function of $V$, and $Y_\Om$ is concave (Theorem \ref{thconpc}), 
the proof trivially follows from Proposition \ref{prop:alternative}. 

\begin{remark}
Though we have succeed in comparing the profiles for small volumes with 
geometric arguments, the global comparison has required global analytic 
properties of the profile.
\end{remark}

\providecommand{\bysame}{\leavevmode\hbox to3em{\hrulefill}\thinspace}
\providecommand{\MR}{\relax\ifhmode\unskip\space\fi MR }
\providecommand{\MRhref}[2]{%
  \href{http://www.ams.org/mathscinet-getitem?mr=#1}{#2}
}
\providecommand{\href}[2]{#2}


\begin{thebibliography}{{\bf GMT}}

\bibitem[{\bf BP}]{BP} 
Christophe~Bavard et Pierre~Pansu, \emph{Sur le volume minimal de
$\rr^{2}$}, Ann. Sci. \'{E}cole. Norm. Sup. \textbf{19} (1986), no.~4, 
479--490. \MR{88b:53048}

\bibitem[{\bf Ba1}]{Bat} 
Vincent~Bayle, \emph{Propri\'et\'es de concavit\'e du profil 
isop\'erim\'etrique et applications}, Th\`ese de Doctorat. 2003.

\bibitem[{\bf Ba2}]{Ba} 
Vincent~Bayle, \emph{A Differential Inequality for the 
Isoperimetric Profile}, Int. Math. Res. Not. (to appear).

\bibitem[{\bf BBG}]{BBG} Pierre B\'erard, G\'erard Besson et Sylvestre 
Gallot, \emph{Sur une in\'egalit\'e isop\'erim\'etrique qui 
g\'en\'eralise celle de Paul L\'evy-Gromov}, Invent. Math. 
\textbf{80} (1985), no.~2, 295-308. \MR{86j:58017}

\bibitem[{\bf BM}]{bemeyer}
Pierre~B{\'e}rard et Daniel~Meyer, \emph{In\'egalit\'es
isop\'erim\'etriques et applications}, Ann.  Sci.  \'{E}cole Norm.  Sup. 
(4) \textbf{15} (1982), no.~3, 513--541.  \MR{84h:58147}

\bibitem[{\bf Bi}]{convexity} 
Richard L.~Bishop, \emph{Infinitesimal convexity implies local convexity}, 
Indiana Univ. Math. J. \textbf{24} (1974-75), 169--172. \MR{MR50:3154}

\bibitem[{\bf Bo}]{bombieri}
Enrico~Bombieri, \emph{Regularity theory for almost minimal currents}, 
Arch. Rational Mech. Anal. \textbf{78} (1982), no.~2, 99--130.  
\MR{MR83i:49077}

\bibitem[{\bf Ch}]{chavel3}
Isaac~Chavel, \emph{Eigenvalues in {R}iemannian geometry}, Pure and
Applied Mathematics, vol.  115, Academic Press Inc., Orlando, FL,
1984.  \MR{86g:58140}

\bibitem[{\bf Ch2}]{chavel} 
\bysame, \emph{Riemannian Geometry: a modern introduction}, 
Cambridge Tracts in Mathematics, no. 108, Cambridge University Press, 
Cambridge, 1993. \MR{MR95j:53001}

\bibitem[{\bf Ch3}]{C2} 
\bysame, \emph{Isoperimetric Inequalities. Differential Geometric and 
Analytic Perspectives}, Cambridge Tracts in Mathematics, no. 145, 
Cambridge University Press, Cambridge, 2001. \MR{2002h:58040}

\bibitem[{\bf Ga}]{G1} 
Sylvestre~Gallot, \emph{In\'egalit\'es isop\'erim\'etriques et
analytiques sur les vari\'et\'es riemanniennes}, Soci\'et\'e 
Math\'ematique de France, Ast\'erisque \textbf{163-164} (1988), 
31-91. \MR{90f:58173}

\bibitem[{\bf Gi}]{giusti}
Enrico~Giusti, \emph{Minimal surfaces and functions of bounded
variation}, Birkh\"auser Verlag, Basel, 1984.  \MR{87a:58041}

\bibitem[{\bf GMT}]{gmt}
Eduardo~Gonzalez, Umberto~Massari, and Italo~Tamanini, \emph{On the
regularity of boundaries of sets minimizing perimeter with a volume
constraint}, Indiana Univ. Math. J. \textbf{32} (1983), no.~1,
25--37.  \MR{84d:49043}

\bibitem[{\bf Gr}]{GR2} 
Misha~Gromov, \emph{Paul L\'evy's Isoperimetric Inequality}, Appendix C 
in Metric Structures for Riemannian and non Riemannian Spaces by M. Gromov, 
Birkhäuser Boston, Inc., Boston, MA, 1999. \MR{2000d:53065}

\bibitem[{\bf G1}]{gruter}
Michael~Gr{\"u}ter, \emph{Boundary regularity for solutions of a
partitioning problem}, Arch.  Rational Mech.  Anal.  \textbf{97}
(1987), no.~3, 261--270.  \MR{87k:49050}

\bibitem[{\bf G2}]{gruter-nuevo}
\bysame, \emph{Optimal regularity for codimension one minimal surfaces
with a free boundary}, Manuscripta Math.  \textbf{58} (1987), no.~3,
295--343.  \MR{88m:49032}

\bibitem[{\bf K}]{K} 
Ernst~Kuwert, \emph{Note on the Isoperimetric Profile 
of a Convex Body}, personal communication.

\bibitem[{\bf M1}]{morgan}
Frank~Morgan, \emph{Geometric measure theory}, third ed., Academic Press
Inc., San Diego, CA, 2000, A beginner's guide.  \MR{2001j:49001}

\bibitem[{\bf M2}]{morgan2}
\bysame, \emph{{Regularity of isoperimetric hypersurfaces in
Riemannian manifolds}}, Trans.  Amer.  Math.  Soc.  (to appear).

\bibitem[{\bf MJ}]{mj}
Frank~Morgan and David~L.~Johnson, \emph{Some sharp isoperimetric
theorems for {R}iemannian manifolds}, Indiana Univ.  Math.  J.
\textbf{49} (2000), no.~3, 1017--1041.  \MR{2002e:53043}

\bibitem[{\bf MR}]{morganritore}
Frank~Morgan and Manuel~Ritor{\'e}, \emph{Isoperimetric regions in
cones}, Trans.  Amer.  Math.  Soc.  \textbf{354} (2002), no.~6,
2327--2339 (electronic).  \MR{2003a:53089}

\bibitem[{\bf RR}]{ritros}
Manuel~Ritor\'{e} and C\'{e}sar~Rosales, \emph{Existence and characterization 
of regions minimizing perimeter under a volume constraint inside Euclidean cones}, 
Trans.~Amer.~Math.~Soc.~(to appear)

\bibitem[{\bf Ro}]{ros}
Antonio~Ros, \emph{The isoperimetric problem}, Lecture series given 
during the \emph{Clay Mathematics Institute Summer School on the Global 
Theory of Minimal Surfaces} at the MSRI, Berkeley, California (2001). 

\bibitem[{\bf RV}]{rosvergasta}
Antonio~Ros and Enaldo~Vergasta, \emph{Stability for hypersurfaces of
constant mean curvature with free boundary}, Geom.  Dedicata
\textbf{56} (1995), no.~1, 19--33.  \MR{96h:53013}

\bibitem[{\bf SZ1}]{zumbrun2bis}
Peter~Sternberg and Kevin~Zumbrun, \emph{Connectivity of phase boundaries 
in strictly convex domains}, Arch. Rational. Mech. Anal., \textbf{141} 
(1998), no.4, 375--400. \MR{99c:49045}  

\bibitem[{\bf SZ2}]{zumbrun2}
\bysame, \emph{On the connectivity of boundaries of sets minimizing
perimeter subject to a volume constraint}, Comm.  Anal.  Geom. 
\textbf{7} (1999), no.~1, 199--220.  \MR{2000d:49062}

\bibitem[{\bf Z}]{bvfunctions}
William~P. Ziemer, \emph{Weakly differentiable functions}, Graduate
Texts in Mathematics, vol.  120, Springer-Verlag, New York, 1989,
Sobolev spaces and functions of bounded variation.  \MR{91e:46046}

\end{thebibliography}
\end{document}